\documentclass[12pt]{article}
\usepackage{mathrsfs, amssymb, eucal, amsmath, amsthm, amscd, }

\makeatletter \@addtoreset{figure}{section}
\def\thefigure{\thesection.\@arabic\c@figure}
\def\fps@figure{h,t}
\@addtoreset{table}{bsection}

\def\thetable{\thesection.\@arabic\c@table}
\def\fps@table{h, t}
\@addtoreset{equation}{section}

\makeatother

\newcommand{\ep}{\varepsilon}
\newcommand{\R}{{\mathbb  R}}

\newcommand{\ds}{\displaystyle}
\newcommand{\ol}{\overline}
\numberwithin{equation}{section}
\newtheorem{thm}{\bf Theorem}[section]

\newtheorem{prop}[thm]{\bf Proposition}

\newtheorem{rem}{\bf Remark}[section]

\newcommand{\p}{\prime}

\newcommand{\la}{\lambda}
\newcommand{\al}{\alpha}

\begin{document}

\newtheorem{theorem}{Theorem}[section]
\newtheorem{definition}[theorem]{Definition}
\newtheorem{lemma}[theorem]{Lemma}
\newtheorem{remark}[theorem]{Remark}
\newtheorem{proposition}[theorem]{Proposition}
\newtheorem{corollary}[theorem]{Corollary}
\newtheorem{example}[theorem]{Example}

\newcommand{\bfi}{\bfseries\itshape}

\newsavebox{\savepar}
\newenvironment{
boxit}{\begin{lrbox}{\savepar}
\begin{minipage}[b]{15.5cm}}{\end{minipage}\end{lrbox}
\fbox{\usebox{\savepar}}}

\title{Periodic orbits in the case of a zero eigenvalue}
\author{Petre Birtea, Mircea Puta, Razvan Micu Tudoran}
\date{}
\maketitle

\begin{abstract}
We will show that if a dynamical system has enough constants of
motion then a Moser-Weinstein type theorem can be applied for
proving the existence of periodic orbits in the case when the
linearized system is degenerate.
\end{abstract}

\section{Introduction.}

Finding periodic solutions of a system of ordinary differential
equations is an old problem in mathematical physics going back to
Lyapunov and Poincare. Periodic solutions were discovered first
for linear conservative systems that appears in mechanics. The
passage from linear to nonlinear systems was taken by Lyapunov
\cite{lyapunov} under the assumption of existence of an integral
of motion and a certain nonresonance condition.

In 1973, Weinstein \cite{weinstein} proved that in the case of a
Hamiltonian system with a positive definite Hamiltonian function
the nonresonance condition is not necessary. Later, Moser
\cite{moser} extended Weinstein's result to the case of a general
dynamical system which posses a constant of motion. More
precisely, let
\begin{equation}\label{system}
\dot x = X(x),
\end{equation}
be a dynamical system generated by the $C^1$ vector field $X$ on a
differentiable manifold $M$ with $x_0$ an equilibrium point, i.e.,
$X(x_0) = 0$. Consider the linearized equations for the
equilibrium point $x_0$,
\begin{equation}
\dot z = DX(x_0)\cdot z.
\end{equation}
Then we have the following result due to Moser \cite{moser}.

\bigskip

\noindent {\bfi Theorem (Moser)} \textit{ Let $I\in C^2$ be an
integral of motion for (\ref{system}) with $dI(x_0)=0$. If
\begin{itemize}
  \item [(i)] $DX(x_0)$ is a non-singular matrix,
  \item [(ii)] $DX(x_0)$ has a pair of pure complex eigenvalues $\pm
  i\omega$ with $\omega \neq 0$,
  \item [(iii)] $d^2I(x_0)$ is positive definite,
\end{itemize} then for sufficiently small $\epsilon$ any integral
surface \[I(x) = I(x_0) + \ep^2\] contains at least one periodic
solution of $X$ whose period is close to the period  of the
corresponding linear system around $x_0$.}

\medskip

The condition $(i)$ of the above theorem implies that the
linearized system around the critical point $x_0$ can not have a
zero eigenvalue. This restriction makes the theorem unapplicable
to a series of examples. We will show that in the case when for
(\ref{system}) one can find enough constants of motion a similar
result can be applied for proving the existence of periodic
orbits. We will also illustrate this with two examples.

\section{The main result.}
\begin{thm}\label{main}
Let $\dot x = X(x)$ be a dynamical system, $x_0$ an equilibrium
point, i.e.,  $X(x_0) = 0$ and $C:= (C_1, \dots, C_k) : M \to
\R^k$ a vector valued constant of motion for the above dynamical
system with $C(x_0)$ a regular value for $C$. If
\begin{itemize}
  \item [(i)] the eigenspace corresponding to the eigenvalue zero of
  the linearized system around $x_0$ has dimension $k$,
  \item [(ii)] $DX(x_0)$ has a pair of pure complex eigenvalues $\pm
  i\omega$ with $\omega \neq 0$,
  \item [(iii)] there exist a constant of motion $I:M\to \R$ for
  the vector field $X$ with $dI(x_0)=0$ and such that
  \[d^2I(x_0)_{|_{W \times W}} >0,\] where
  $W=\bigcap\limits_{i=1}^k \ker d C_i (x_0),$
\end{itemize} then for each
  sufficiently small $\ep \in \R$, any integral surface \[I(x) =
  I(x_0) + \ep^2\] contains at least one periodic solution of $X$
  whose period is close to the period  of the corresponding linear
  system around $x_0$.
\end{thm}

\begin{proof}

If $C_i \in C^\infty (M, \R)$ is a constant of motion for the
dynamic generated by the vector field $X$ then $DX(x_0)\nabla C_i
(x_0)=0$, and hence $\nabla C_i (x_0) \in \ker DX(x_0)$.

Because $C(x_0)$ is a regular value for $C$ we have that
$dC_i(x_0)$, $i=\ol{1, k}$ are linearly independent vectors in the
tangent space $T_{x_0}M$. Then, hypothesis $(i)$ and the fact that
$C_1, \dots, C_k \in C^\infty (M, \R)$ are constants of motion for
$X$ implies the following equality,
\[span \{\nabla C_i (x_0) : i = \ol{1,k}\}=\ker DX (x_0)
(=V_{\la =0})\] where $V_{\la=0}$ is the eigenspace corresponding
to the zero eigenvalue of the matrix which is canonically
associated to the linear part at the equilibrium of interest $x_0$
of our system determined by $X$.

This argument implies that the reduced system
\begin{equation*}
\left \{
\begin{array}{l}
\dot{x}= X(x) \\
C(x)=C(x_0)
\end{array}
\right.
\end{equation*} which is the original system
restricted to the submanifold $C^{-1}(C(x_0))$ has the linearization
about $x_0$ without eigenvalue zero.

The function $I_{|C^{-1}(C(x_0))}:C^{-1}(C(x_0)) \rightarrow \R$ is
a first integral for the reduced system with $d(I_{|C^{-1}(C(x_0))})
(x_0) =0$ and hypothesis $(iii)$ obviously implies that
$d^2(I_{|C^{-1}(C(x_0))})(x_0)>0$. By the Moser theorem we have that
for sufficiently small $\ep \in \R$, any integral surface
\[I(x) = I(x_0) + \ep^2\] contains at least one periodic solution
of the reduced system and hence of the initial system.
\end{proof}

\medskip

\begin{rem}
If the dynamic (\ref{system}) is Hamilton-Poisson and $x_0$ is
regular point in the sense that it is contained in a maximal
dimension symplectic leaf of $(M, \{ \})$ which is determined by
the Casimirs $C_1, \dots, C_k$, then by the theorem of Weinstein
\cite{weinstein} one has the existence of  $(dim P-k)/2$ periodic
orbits.
\end{rem}

\section{Examples.}

\noindent {\bfi Rigid body with one control.} Let us consider the
rigid body dynamics with one control,
\begin{equation}\label{e1}
    \left\{\begin{array}{l} \dot m_1 = a_1 m_2 m_3\\ \dot m_2 = a_2
    m_1 m_3\\ \dot m_3 = (a_3 - l) m_1 m_2\end{array}\right.
\end{equation} where $l \in \R$ is the gain parameter.

Let us make now the following notation $\alpha:=
\ds\frac{a_3-l}{a_3}$. Then it is not hard to see that our
dynamics \eqref{e1} has the following Hamilton-Poisson realization
$(\R^3, \Pi_\alpha, H_\alpha)$, where
\[\Pi_\al \stackrel{def}{=}\left[\begin{array}{rrr}0&-m_3& \al
m_2\\m_3& 0 & - \al m_1\\ -\al m_2 & \al m_1 & 0
\end{array}\right]
\] is the Poisson structure and $H_\al (m_1, m_2, m_3) \stackrel{def}{=} \ds\frac{1}{2}
\left(\ds\frac{m_1^2}{I_1}+\ds\frac{m_2^2}{I_2}+\ds\frac{m_3^2}{\al
I_3}\right)$ is the Hamiltonian function. Moreover, the smooth
function $C_\al \in C^\infty (\R^3, \R)$ given by \[C_\al (m_1, m_2,
m_3) \stackrel{def}{=} \al m_1^2 + \al m_2^2 + m_3^2\] is a Casimir
of our Poisson configuration $(\R^3, \Pi_\al)$.

Let us concentrate now to the equilibrium state \[e_1^M = (M,0,0),
\; M \in \R^*\] of our dynamics \eqref{e1}. Then under the
restriction $l < a_3$ we have succesivelly,
\begin{itemize}
  \item [(i)] The restriction of the dynamics \eqref{e1} to the
  coadjoint orbit
  \begin{equation}\label{e2}
    \al m_1^2 + \al m_2^2 + m_3^2 = \al M^2
  \end{equation} gives rise to a Hamiltonian system
  on a symplectic manifold.

  \item [(ii)] $span \ (\nabla C_\al (e_1^M)) = V_{\la=0} = span
  \left(\left[\begin{array}{c} 1\\0\\0\end{array}\right]\right)$
  where \[V_{\la=0} =\left\{\left.\left[\begin{array}{c}
  m_1\\m_2\\m_3\end{array}\right]\in \R^3 \ \right| A(e_1^M)
  \left[\begin{array}{c} m_1\\m_2\\m_3\end{array}\right] =
  \left[\begin{array}{c}0\\0\\0\end{array}\right] \right\},\]
  $A(e_1^M)$ being the matrix of the linear part of the dynamics
  \eqref{e1} at the equilibrium of interest $e_1^M$, $M \in \R^*$.

  \item [(iii)] The matrix of the linear part of our reduced
  dynamics to \eqref{e2} has at the equilibrium $e_1^M$ the
  following characteristic roots:
  \[\la_{1,2} = \pm Mi \sqrt{-a_2 (a_3 - l)}\].

  \item[(iv)] The smooth function $F_{\frac{1}{\al I_1}} \in
  C^\infty (\R^3, \R)$ given by: \[F_{\frac{1}{\al I_1}} (m_1, m_2,
  m_3) = \ds\frac{1}{2}
  \left(\ds\frac{m_1^2}{I_1}+\ds\frac{m_2^2}{I_2}+\ds\frac{m_3^2}{\al
  I_3}\right) - \ds\frac{1}{2\al I_1} (\al m_1^2 + \al m_2^2 +
  m_3^2)\] is a constant of motion and $e_1^M$ is a local minimum of
  $F_{\frac{1}{al I_1}}$ with the constraint \eqref{e2}.
\end{itemize}

 Then via Theorem \ref{main} we have:

\begin{prop}
If $l < a_3$ then the reduced dynamics to the coadjoint orbit
\eqref{e2} has near the equilibrium  state $e_1^M$, $M \in \R^*$
at least one periodic solution whose period is close to
\[ \ds\frac{2 \pi}{\mid M \mid \sqrt{-a_2 (a_3 - l)}}.\]\hfill{$\Box$}
\end{prop}

\begin{rem}
 Similar results can  be also obtained for the equilibrium states
 \[e_2^M = (0,M,0), \; M \in \R^*\] and \[e_3^M = (0,0,M), \; M \in
 \R^*.\] \hfill{$\Box$}
\end{rem}

\bigskip

\noindent {\bfi Clebsch system.} It is well known that the Clebsch
system can be written in the following form:
\begin{equation}\label{e3}
    \left\{\begin{array}{l} \dot x_1 = x_2 p_3 - x_3 p_2\\ \dot x_2
    = x_3 p_1 - x_1 p_3\\ \dot x_3 = x_1 p_2 - x_2 p_1\\ \dot p_1 =
    (a_3 - a_2) x_2 x_3\\ \dot p_2 = (a_1 - a_3) x_1 x_3\\ \dot p_3
    = (a_2 - a_1) x_1 x_2 \end{array}\right.
\end{equation} where \[\begin{array}{l} a_1, \ a_2,\ a_3 \in \R\\
a_1>0, \ a_2 >0, \ a_3 >0\\a_1 \not = a_2 \not = a_3\end{array}\]
(see for details Dubrovin, Krichever and Novikov \cite{dubrovin}).

It is not hard to see that the smooth functions $H, C, D \in
C^\infty (\R^6, \R)$ given by:
\begin{eqnarray*}
H(x_1, x_2, x_3, p_1, p_2, p_3) &=& \ds\frac{1}{2}
(a_1x_1^2+a_2x_2^2+a_3x_3^2+p_1^2+p_2^2+p_3^2)\\ C(x_1, x_2, x_3,
p_1, p_2, p_3) &=& \ds\frac{1}{2} (x_1^2+x_2^2+x_3^2)\\D(x_1, x_2,
x_3, p_1, p_2, p_3) &=& x_1p_1+x_2p_2+x_3p_3\end{eqnarray*} are
constants of motion for the Clebsch system.

Let us concentrate now to the equilibrium state $e_1^M =
(M,0,0,0,0,0)$, $M \in \R^*$. Then under the restrictions:
\begin{equation*}
a_3>a_1 \text{ and } a_2 > a_1
\end{equation*}
we have successively,
\begin{itemize}
  \item [(i)] $span \ (\nabla C_\al (e_1^M), \nabla D(e_1^M)) = V_{\la=0} $
  where \[V_{\la=0} =\left\{\left.\left[\begin{array}{c}
  m_1\\m_2\\m_3\\p_1\\p_2\\p_3\end{array}\right]\in \R^6 \ \right| A(e_1^M)
  \left[\begin{array}{c} m_1\\m_2\\m_3\\p_1\\p_2\\p_3\end{array}\right] =
  \left[\begin{array}{c}0\\0\\0\\0\\0\\0\end{array}\right] \right\},\]
  $A(e_1^M)$ being the matrix of the linear part of the dynamics
  \eqref{e3} at the equilibrium $e_1^M$.

  \item [(ii)] The matrix of the linear part of our reduced
  dynamics to the constraint
  \begin{equation}\label{e4}
  \left\{ (x_1, x_2, x_3, p_1, p_2, p_3) \in \R^6 \ \left|
  \begin{array}{l} x_1^2 + x_2^2 + x_3^2 = M^2\\ x_1p_1 + x_2p_2 +
  x_3 p_3=0\end{array}\right.\right\},
  \end{equation} at the equilibrium $e_1^M$ has the following
  characteristic roots:
  \[\begin{array}{l} \la_{1,2} = \pm iM \sqrt{a_3  - a_1},\\ \la_{3,4}
  = \pm iM \sqrt{a_2 - a_1}.\end{array}\]

  \item[(iii)] The smooth function $F_{a_1} \in
  C^\infty (\R^6, \R)$ given by: \begin{eqnarray*} F_{a_1} (x_1, x_2,
  x_3, p_1, p_2, p_3) &=& \ds\frac{1}{2}
  \left(a_1x_1^2+a_2x_2^2+a_3x_3^2 +p_1^2 + p_2^2 + p_3^2\right)\\& -
  &\ds\frac{a_1}{2} (x_1^2 + x_2^2 +
  x_3^2)\end{eqnarray*}
   is a constant of motion and $e_1^M$ is a local minimum of
  $F_{a_1}$ with the constraint \eqref{e4}.
\end{itemize}

Then via Theorem \ref{main} we have:

\begin{prop}
If $a_2 < a_1$ and $a_3 > a_1$ then the reduced dynamics to
\eqref{e4} has near  $e_1^M$, $M \in \R^*$ at least one periodic
solution whose period is close to $ \ds\frac{2 \pi}{\mid M\mid
\sqrt{a_3 - a_1}}$ and $\ds\frac{2\pi}{\mid M \mid
\sqrt{a_2-a_1}}$.\hfill{$\Box$}
\end{prop}

\begin{rem}
 Similar results can  be also obtained for the equilibrium states:
 \[e_2^M = (0,M,0,0,0,0), \; M \in \R^*\] and \[e_3^M = (0,0,M,0,0,0), \; M \in
 \R^*.\] \hfill{$\Box$}
\end{rem}

\end{document}